\documentclass[12pt]{article}
\usepackage{graphicx,psfrag,epsfig, color,float}
\usepackage{amssymb,amsmath,amscd,amsthm}
\usepackage{graphicx,psfrag,epsfig}

\usepackage{graphicx}
\usepackage[active]{srcltx}

\newtheorem{theorem}{Theorem}[section]

\newtheorem{lemma}[theorem]{Lemma}

\date{}

\begin{document}

\date{}
\title{Metastable distributions of Markov chains with rare transitions}
\author{
 M. Freidlin\footnote{Dept of Mathematics, University of Maryland,
College Park, MD 20742, mif@math.umd.edu}, L.
Koralov\footnote{Dept of Mathematics, University of Maryland,
College Park, MD 20742, koralov@math.umd.edu}
} \maketitle

\begin{abstract}
In this paper we consider Markov chains with transition rates that depend on a small parameter $\varepsilon$.
Under a mild assumption on the asymptotics of these transition rates, we describe the behavior of the chain at various $\varepsilon$-dependent time scales, i.e., we calculate the metastable distributions.
\end{abstract}

{2010 Mathematics Subject Classification Numbers: 37A25, 37A30, 60J27,

60J28, 60F10.}
\\


\section{Introduction.} \label{intro}

Consider a family $X^\varepsilon_t$ of Markov chains on a state space $S = \{1,...,N\}$, where $\varepsilon$ is a small parameter. The time may be
continuous or discrete -- we start by considering the case when $t \in \mathbb{R}^+$. The case of discrete time is similar  and is briefly discussed in Section~\ref{gene}. Let $q_{ij}(\varepsilon)$, $i, j \in S$, $i \neq j$,  be the transition rates, i.e.,
\begin{equation} \label{contti}
\mathrm{P}(X^\varepsilon_{t + \Delta} = j | X^\varepsilon_t = i) = q_{ij}(\varepsilon) (\Delta + o(\Delta))~~~{\rm as}~\Delta \downarrow 0,~~i \neq j.
\end{equation}
We will be interested in the behavior of $X^\varepsilon_t$ as $\varepsilon \downarrow 0$ and, simultaneously, $t = t(\varepsilon) \rightarrow~\infty$.
The results on the asymptotic behavior of $X^\varepsilon_t$ can be viewed as a refinement of the ergodic theorem for Markov chains (which concerns the asymptotics with respect to the time variable only) and are obviously closely related to the spectral properties of the transition matrix.
The double limit at hand depends on how the point $(1/\varepsilon, t(\varepsilon))$ approaches infinity. Roughly speaking, one can divide the neighborhood of infinity into a finite number of domains such that $X^\varepsilon_{t(\varepsilon)}$ has a limiting distribution (which depends on the initial point) when $(1/\varepsilon, t(\varepsilon))$ approaches infinity
without leaving a given domain. For different domains, these limits are different. These will be referred to as metastable distributions.

Families of parameter-dependent Markov chains arise in a variety of applications. In particular, this is a natural object in problems concerning random perturbations of dynamical systems (\cite{FW}). If a dynamical system has $N$ asymptotically stable attractors, each attractor (or rather its small neighborhood) can be associated with a state of a Markov chain, while the transition times between different states are due to large deviations and are determined by the action functional of the perturbed system. In this example, these transition times are exponentially large, i.e., the transition rates are exponentially small with respect to the size of the perturbation~$\varepsilon$.

As shown in \cite{F5}, \cite{FW}, the long-time behavior of the perturbed process can typically be understood using the notion of the hierarchy of cycles.  The hierarchy of cycles means, roughly speaking, that for each $0 \leq r \leq \rho$
the set of attractors is decomposed into disjoint subsets $C^{r}_1,...,C^{r}_{n_r}$ (cycles of rank $r$), up to the maximal rank $\rho < N$. The individual attractors are the cycles of rank zero, they are combined in disjoint sets - cycles of rank one, those are combined in cycles
of rank two, etc., until the cycle of maximal rank $\rho < N$ containing all the attractors. With probability close to one as $\varepsilon
\downarrow 0$, the process goes from a neighborhood of one of the attractors in $C^{r}_i$ to a small neighborhood of one the
attractors of the next cycle $C^{r}_j$, thus remaining within a cycle of rank $r+1$. The transition time between $C^{r}_i$ and $C^{r}_j$
is determined by the asymptotics, as $\varepsilon \downarrow 0$,  of the transition rates between individual attractors that belong to the union of these cycles. The process leaves the cycle of rank $r+1$ only after an exponentially large number of transitions between cycles of rank $r$.

For each  $\lambda$ (except a finite number of values) and the time scale $t(\varepsilon) \sim \exp(\lambda/\varepsilon)$, with probability close to one, the process can be found
in a neighborhood of a particular attractor (meta-stable state), which depends on the initial state. The meta-stable state is a piece-wise constant function of the parameter~$\lambda$.

The above description with the hierarchy of cycles and the meta-stable states is valid, however, only if the notion of the unique ``next" cycle (and, consequently, unique meta-stable state)  can be correctly defined, which is not the case in many interesting situations. For example, when the unperturbed dynamics has certain symmetries 
(or ``rough symmetries" as in \cite{FF}), the notion of ``next" is not defined uniquely even for individual attractors (cycles of rank zero). A similar phenomenon was observed in  \cite{FKW}, due not to symmetries but to degeneration of the unperturbed dynamics in a part of the phase space. Other systems leading to parameter-dependent Markov chains with no unique metastable state arise in the study of  various models of non-equilibrium statistical mechanics at low temperatures
(see \cite{MNOS} and references therein). The transition rates for such chains typically decay exponentially in the parameter corresponding to inverse temperature. Various large deviation results for Markov chains with exponentially small transition rates were obtained in \cite{WE} \cite{OS1}, \cite{OS2}. 

Yet another example is provided by dynamical systems with heteroclinic networks. Namely, assume that there are finitely many stationary points with heteroclinic connections that together form a connected set. Assume that the entire network serves as an attractor for the dynamical system. The flow lines (or their sufficiently small neighborhoods) connecting pairs of stationary points can be associated with the states of a Markov chain. After a random perturbation of size $\varepsilon$, the transition times between the neighborhoods of such flow lines scale as powers of $\varepsilon$ (up to a logarithmic factor). So the notions of
the hierarchy of cycles and the meta-stable states could apply (at times that scale as powers of $\varepsilon$, rather than exponentially). However,
the notion of ``next" state may again be not defined uniquely, since the exit from a neighborhood of a heteroclinic connection (say, between stationary points $A$ and $B$) can happen along either of the heteroclinic connections leading out of $B$. A detailed study of motion along heteroclinic networks is the subject of \cite{BP}.  It should be stressed that the dynamics in this example, as well as in the case of asymptotically stable attractors discussed above, is only approximately described by a Markov chain of the type considered in the current paper. A reduction of the true dynamics to a finite-state Markov chain requires non-trivial analysis.

In the current paper, we introduce the notion of the hierarchy of Markov chains in a general setting. We show that it should replace the notion of the hierarchy of cycles. The meta-stable states are replaced by meta-stable distributions. We do not require the transition rates between different states to scale exponentially (or have any specific asymptotic behavior), but only assume that there is a certain asymptotic relation between the ratios of transition rates.

More precisely, we will say that the family $X^\varepsilon_t$ is {\it asymptotically regular} as $\varepsilon \downarrow 0$ if the following two conditions hold.

(a) The transition rates $q_{ij}(\varepsilon)$ are positive\footnote{The main result can be obtained even if the  positivity assumption is relaxed, as mentioned in Section~\ref{gene}.} for $\varepsilon > 0$ and all $i \neq j$.

(b) For each $i, j, k, l\in S$ satisfying $i \neq j$, $k \neq l$, the following finite or infinite limit exists
\[
\lim_{\varepsilon \downarrow 0} ({q_{ij}(\varepsilon)}/{q_{kl}(\varepsilon)})  \in [0, \infty].
\]
%

In Section~\ref{hie}, we will introduce the hierarchy of Markov chains. The hierarchy will be defined inductively by successively reducing the state space, i.e., combining the elements of the state space into subsets that will serve as states for the chain of higher rank. In order to perform an inductive step, we will require each chain appearing in the construction to be asymptotically regular. While this condition may seem not quite explicit,  we will show that it holds if the transition rates of the original family satisfy a relatively simple condition.

In Section~\ref{are}, we discuss some simple properties of asymptotically regular Markov chains.
In Section~\ref{se3}, we formulate and prove the main theorem on the meta-stable behavior of the original process.
In Section~\ref{car}, we prove that a condition on the transition rates of the family of Markov chains guarantees that all the chains in the hierarchy
are asymptotically regular.
We briefly discuss a couple of generalizations in Section~\ref{gene}. 

\section{Hierarchy of Markov chains. } \label{hie}

\subsection{Reduced Markov chain.} \label{rmc}

Given an asymptotically regular family of Markov chains $X^\varepsilon_t$ with $N \geq 2$, we will construct a reduced Markov chain (later also referred to as the
reduced Markov chain of rank one).
First, we define a discrete-time  Markov chain on $S$, which will be referred to as the {\it skeleton  Markov chain} and denoted by $Z_n$. Its
transition probabilities are defined by
\[
P_{ij} = \lim_{\varepsilon \downarrow 0} (q_{ij}(\varepsilon)/\sum_{j' \neq i} q_{ij'}(\varepsilon)),~~j \neq i;~~~P_{ii} = 0.
\]
Observe that the above limit exists since the family of chains $X^\varepsilon_t$  is asymptotically regular.
Recall (see \cite{Doob}) that the (finite) state space of a Markov chain can be uniquely decomposed into a disjoint union of ergodic classes and the sets consisting of individual transient states,
\begin{equation} \label{decomp}
S = S_1 \bigcup ... \bigcup S_n.
\end{equation}
Note that $n < N$ since each ergodic class of $Z_n$ contains at least two states, which follows from the fact that for each $i$ there is $j \neq i$ such that $P_{ij} > 0$.

Next, for each $1 \leq k \leq n$, we define the Markov chains $Y^{k, \varepsilon}_t$ by narrowing the state space $S$ to $S_k$ and retaining the same
transition rates for $i, j \in S_k$ as in the original Markov chain $X^\varepsilon_t$. Let $\mu^k(i,\varepsilon)$ be the invariant measure of the state $i \in S_k$ for the chain~$Y^{k, \varepsilon}_t$.

Finally, we define the reduced Markov chain. Its state space is the set $\{S_1,..., S_{n}\}$.
The transition rate between $S_k$ and $S_l$, $k \neq l$, denoted by $Q_{kl}(\varepsilon)$, is defined by
\begin{equation} \label{defq0}
Q_{kl}(\varepsilon) = \sum_{i \in S_k} \sum_{j \in S_l} \mu^{k}(i, \varepsilon) q_{ij}(\varepsilon).
\end{equation}

\subsection{Definition of the hierarchy.}

We will use induction to define the  {\it reduced Markov chains} $X^{r, \varepsilon}_t$, $0 \leq r \leq \rho $, with some $0 \leq \rho < N$. The reduced Markov chain of rank zero (i.e., corresponding to $r = 0$) will coincide with $X^\varepsilon_t$, while the reduced Markov chain of rank $\rho$ will be trivial (i.e., its state space will contain one element).
The entire collection of reduced Markov chains will be referred to as the hierarchy.

For each $0 \leq r \leq \rho-1$, by partitioning the state space of $X^{r, \varepsilon}_t$ into $n_{r+1}$ subsets (referred to as clusters), we will define Markov chains $Y^{r, k, \varepsilon}_t$, $1 \leq k \leq n_{r+1}$.
The chain $Y^{r, k, \varepsilon}_t$ will be referred to as the $k$-th chain of rank $r$.

We set  $n_0 = N$ and  $S^0_1 = 1$,..., $S^0_{n_0} = n_0$. These are clusters of rank zero. The reduced Markov chain of rank zero, denoted by $X^{0, \varepsilon}_t$, is defined on the state
space $S^0 = \{S^0_1,..., S^0_{n_0}\}$ and has transition rates $Q^{0}_{ij}(\varepsilon) = q_{ij}(\varepsilon)$. Thus it coincides with the original Markov chain.

If $N = 1$, this results in the trivial hierarchy. If $N \geq 2$, we can apply the above construction of the reduced Markov chain. We set $n_1 = n$ and use the following notation
\[
S = S^1_1 \bigcup ... \bigcup S^1_{n_1}.
\]
for the decomposition of the original state space into ergodic classes and sets containing individual transient states for the skeleton chain. The sets $S^1_1,...,S^1_{n_1}$ will be referred to as clusters of rank one. For each $1 \leq k \leq n_1$, we set $Y^{0, k, \varepsilon}_t = Y^{k, \varepsilon}_t$, which will be be referred to as the $k$-th Markov chain of rank zero. Let $\mu^{0,k}(i, \varepsilon)$ be the invariant measure of the state $S^{0}_i \in S^1_k$ for the chain $Y^{0, k, \varepsilon}_t$ defined above.

The reduced Markov chain introduced above will also be referred to as the reduced Markov chain or rank one and denoted by $X^{1, \varepsilon}_t$. Its state space is the set $\{S^1_1,..., S^1_{n_1}\}$. The transition rate between $S^1_k$ and $S^1_l$, $k \neq l$, denoted by $Q^1_{kl}(\varepsilon)$,  is defined, in conformance with (\ref{defq0}), by
\begin{equation} \label{defq}
Q^1_{kl}(\varepsilon) = \sum_{i \in S^1_k} \sum_{j \in S^1_l} \mu^{0,k}(i, \varepsilon) Q^0_{ij}(\varepsilon).
\end{equation}
If $X^{1, \varepsilon}_t$ is asymptotically regular, we can
replicate the above construction, i.e., define the skeleton chain corresponding to $X^{1, \varepsilon}_t$ and partition $\{S^1_1,...,S^1_{n_1} \}$ into clusters $S^2_1,...,S^2_{n_2}$. For each $1 \leq k \leq n_2$, we can define the $k$-th Markov chain of rank one, denoted by $Y^{1, k, \varepsilon}_t$. The reduced Markov chain $X^{2, \varepsilon}_t$ of rank two is defined on the state space $\{S^2_1,...,S^2_{n_2}\}$. The construction then continues inductively, assuming that all the reduced chains are asymptotically regular.

\vskip -5pt
\begin{figure}[htbp]
    \centerline{\includegraphics[height=3.6in, width= 5.5in,angle=0]{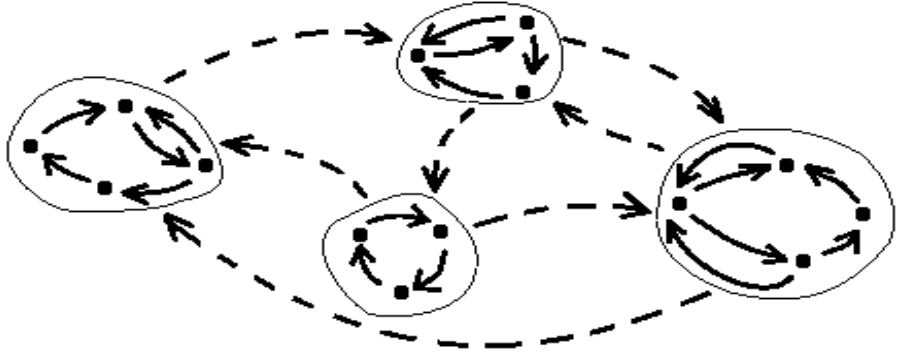}}
\vskip -70pt
  \caption{In this example, $n_0 = N = 14$, $n_1 = 4$, $n_2 =1$, and $\rho =2$.  The solid arrows denote non-zero transitions for the skeleton Markov chain corresponding to $X^{0,\varepsilon}_t$. The dashed arrows denote non-zero transitions for the skeleton Markov chain corresponding to
$X^{1,\varepsilon}_t$.}
    \label{secondkind}
\end{figure}

Let $\rho$ be the first index such that $n_{\rho} = 1$. Observe that $\rho < N$ since $n_0 = N$ and $n_{r+1} < n_r$ for each $r$. Since $n_{\rho} = 1$, there is only one Markov chain of rank $\rho-1$, and the reduced Markov chain of rank $\rho$ is trivial -- its state space consists of one element $S^\rho_1$.
This completes the construction of the hierarchy.
\\

The reduced Markov chain of rank $r$, $0 \leq r \leq \rho$,  will be denoted by $X^{r,\varepsilon}_t$, its transition rates will be denoted by $Q^r_{kl}(\varepsilon)$, the $k$-th Markov chain of rank $r$, $0 \leq r < \rho$, $1 \leq k \leq n_{r+1}$, will be denoted by $Y^{r, k, \varepsilon}_t$, and the invariant measure of a state $S^r_i \in S^{r+1}_k$ for this Markov chain will be denoted by $\mu^{r, k}(i, \varepsilon)$.

For a state $j \in S$, we'll write that $j \prec S^{r}_k$ if there is a sequence $j = j_0, j_1,...,j_r = k$ such that
\begin{equation} \label{inclu}
j = S^0_{j_0} \in S^1_{j_1} ... \in S^r_{j_r}.
\end{equation}
For $1 \leq k \leq n_r$ and $j$ such that $j \prec S^r_k$ does not hold, we will need the functions $\widetilde{Q}^r_{kj}(\varepsilon)$. These are defined inductively, namely,
$\widetilde{Q}^0_{kj}(\varepsilon) = {Q}^0_{ij}(\varepsilon)$ if $S^0_k = i$ and
\[
\widetilde{Q}^r_{kj}(\varepsilon) = \sum_{i: S^{r-1}_i \in S^r_k} \mu^{r-1,k}(i, \varepsilon) \widetilde{Q}^{r-1}_{ij}(\varepsilon),~~1 \leq r \leq \rho-1.
\]
Intuitively, these serve as transition rates from a cluster of rank $r$ to an individual state, although, unlike ${Q}^r_{kl}(\varepsilon)$, they
don't correspond to transition rates of any of the Markov chains introduced above.
\\

Let us stress that the inductive construction of the hierarchy is possible under the condition that all the reduced chains that appear at each step are asymptotically
regular.   We will say that $X^\varepsilon_t$ is {\it completely asymptotically regular} if
for each $a > 0$ and each $(i_1,...,i_a)$,
$(j_1,...,j_a)$, $(k_1,...,k_a)$, and $(l_1,...,l_a)$  the following finite or infinite limit exists
\begin{equation} \label{care}
\lim_{\varepsilon \downarrow 0} \left( \frac{q_{i_1 j_1}(\varepsilon)}{q_{k_1 l_1}(\varepsilon)} \times ... \times \frac{q_{i_a j_a}(\varepsilon)}{q_{k_a l_a}(\varepsilon)}\right) \in [0, \infty],
\end{equation}
provided that $i_1 \neq j_1,...,i_a \neq j_a, k_1 \neq l_1,...,k_a \neq l_a$.

\begin{lemma} \label{lereg0} Suppose that $X^\varepsilon_t$ is completely asymptotically regular.
Then the reduced chain is also completely asymptotically regular.
\end{lemma}
This lemma will be proved in Section~\ref{car}. For now, we observe that  if $X^\varepsilon_t$ is completely asymptotically regular,
then, by Lemma~\ref{lereg0}, so are the reduced Markov chains that appear at each step of the inductive construction of the hierarchy, which implies that all of them are asymptotically regular.
\\

\section{Asymptotically regular families of Markov chains.} \label{are}

In order to prepare for the discussion of metastability, we need several lemmas on asymptotically regular families of Markov chains. Let $\mu(i, \varepsilon)$ denote the invariant measure of the state $i$ in an asymptotically regular Markov chain $X^\varepsilon_t$ with $N \geq 2$. We'll be interested in the asymptotics of $\mu(i, \varepsilon)$ in the case when the skeleton chain
has one ergodic class and no transient states (in which case $\rho = 1$ and no additional assumptions are required in order to define the hierarchy).
Note that this condition is satisfied for each of the Markov chains $Y^{r, k, \varepsilon}_t$ defined above, with the exception that the number of
states for such a chain may be equal to one.

Let $\lambda$ be the invariant measure of the skeleton chain $Z_n$.
Define
\[
T(i, \varepsilon) = ({\sum_{j \in S, j \neq i} q_{ij}(\varepsilon)})^{-1},~~~\bar{T}(\varepsilon) = \sum_{i' \in S}( \lambda(i') T(i',\varepsilon)).
\]
Thus $T(i,\varepsilon)$ is the average of the exponentially distributed exit time of $X^\varepsilon_t$ from the state~$i$. The function $\bar{T}(\varepsilon)$ is  the average time it takes $X^\varepsilon_t$ to make one step,
where the average is calculated with respect to the invariant measure for the skeleton chain.

\begin{lemma} \label{hjjh} Let $X^\varepsilon_t$ be asymptotically regular and $N \geq 2$. Suppose that the skeleton chain has one ergodic class and no transient states.  Then
\[
\mu(i, \varepsilon) \sim \lambda(i) T(i, \varepsilon)/ \bar{T}(\varepsilon)~~~{\it as}~\varepsilon \downarrow 0.
\]
\end{lemma}
\proof Let $t^\varepsilon_1, t^\varepsilon_2,...$ be the times of the jumps of the process $X^\varepsilon_t$. Let $Z^\varepsilon_n$ be the discrete-time Markov chain obtained from $X^\varepsilon_t$ by discretizing time, i.e., $Z^\varepsilon_0 = X^\varepsilon_0$, $Z^\varepsilon_n = X^\varepsilon_{t_n}$, $n \geq 1$.  Let $P^\varepsilon_{ij}$ be the transition probabilities for $Z^\varepsilon_n$ and $\lambda^\varepsilon$ be its invariant measure. Then, from the definition of the skeleton chain, it follows that
\begin{equation} \label{prox}
\lim_{\varepsilon \downarrow 0} P^\varepsilon_{ij} = P_{ij},~~i \neq j,
\end{equation}
where $P_{ij}$ are the transition probabilities for $Z_n$. Since $Z_n$ has one ergodic class and no transient states, this implies that
\[
\lim_{\varepsilon \downarrow 0} \lambda^\varepsilon(i)  = \lambda(i) > 0,~~ i \in S.
\]
By the Law of Large Numbers for Markov chains, $\mu(i, \varepsilon)$ is equal to the asymptotic (as $t \rightarrow \infty$) proportion of time that the process $X^\varepsilon_t$ spends in the state $i$. Therefore,
\[
\mu(i, \varepsilon) = \lambda^\varepsilon(i)  T(i, \varepsilon)/ \sum_{i' \in S}( \lambda^\varepsilon(i') T(i',\varepsilon)).
\]
Combining the latter two equalities, we obtain the result claimed in the lemma.
\qed
\\

From the asymptotic regularity of $X^\varepsilon_t$ it follows that there is a limiting probability measure
\begin{equation} \label{invme}
\mu(i) = \lim_{\varepsilon \downarrow 0} \mu(i, \varepsilon) \in [0,1],~~i \in S.
\end{equation}

Given two functions $t(\varepsilon), s(\varepsilon): (0,\infty) \rightarrow (0,\infty)$, we'll write $s(\varepsilon) \ll t(\varepsilon)$ if  $s(\varepsilon) = o(t(\varepsilon))$ as $\varepsilon \downarrow 0$. Let $\alpha(j, t)$ be the proportion of time, prior to $t$, that the process spends in $j$.
\begin{lemma} \label{interm} Let $X^\varepsilon_t$ be asymptotically regular and $N \geq 2$. Suppose that the skeleton chain has one ergodic class and no transient states. Suppose that $t(\varepsilon)$ is such that $\bar{T}(\varepsilon) \ll t(\varepsilon)$. Then for each $i, j \in S$,
\begin{equation} \label{w1}
\lim_{\varepsilon \downarrow 0} {\mathrm{P}_i}({X}^\varepsilon_{t(\varepsilon)} = j ) =  \mu(j).
\end{equation}
\begin{equation} \label{w2}
{\mathrm{E}_i} \alpha(j, t(\varepsilon)) \sim  \mu(j, \varepsilon),~~~{\it as}~\varepsilon \downarrow 0.
\end{equation}
For each $c > 0$, there are $\delta(c) > 0$ and $\varepsilon_0 > 0$ such that
\begin{equation} \label{w3}
{\mathrm{P}_i}(\alpha(j, t(\varepsilon)) < (1-c)\mu(j, \varepsilon) ) \leq e^{-\delta(c) t(\varepsilon)/\bar{T}(\varepsilon)},~~~\varepsilon \leq \varepsilon_0.
\end{equation}
\begin{equation} \label{w3b}
{\mathrm{P}_i}(\alpha(j, t(\varepsilon)) > (1+c)\mu(j, \varepsilon) ) \leq e^{-\delta(c) t(\varepsilon)/\bar{T}(\varepsilon)},~~~\varepsilon \leq \varepsilon_0.
\end{equation}
Here the subscript $i$ stands for the initial location of the process.\footnote{Formula (\ref{w1}) can be improved to ${\mathrm{P}_i}({X}^\varepsilon_{t(\varepsilon)} = j ) \sim  \mu(j, \varepsilon)~{\rm as}~\varepsilon \downarrow 0$, but we don't need it here.}
\end{lemma}
\proof Let $i^* \in S$ be such that $\mu(i^*) > 0$. Given $i \in S$, find a sequence $i_0, i_1,..., i_k$ such that
$i_0 = i$, $i_k = i^*$, $0 \leq k < N$, and $P_{i_0 i_1},...,P_{i_{k-1} i_k} > 0$, where the
latter are the transition probabilities for the skeleton Markov chain. Then, examining the transition rates of $X^\varepsilon_t$, it is easy to see that
\[
\mathrm{P}_i(Z^\varepsilon_0 = i_0,...,Z^\varepsilon_k = i_k,~~t^\varepsilon_k < \bar{T}(\varepsilon) < t^\varepsilon_{k+1}) \geq a
\]
for some positive constant $a$ and all sufficiently small $\varepsilon$. In particular, $\mathrm{P}_i(X^\varepsilon_{\bar{T}(\varepsilon)} = i^*) \geq a$. Since $i$ was arbitrary, this provides an upper bound on the speed of convergence of $X^\varepsilon_t$ to the invariant distribution, i.e., ${\mathrm{P}_i}({X}^\varepsilon_{t(\varepsilon)} = j ) - \mu(j, \varepsilon) \rightarrow 0$ as $\varepsilon \downarrow 0$ if $\bar{T}(\varepsilon) \ll t(\varepsilon)$. Combined with (\ref{invme}), this implies~(\ref{w1}).

Let $n(\varepsilon) = t(\varepsilon)/\bar{T}(\varepsilon)$. Let $\beta(j, t)$ be the amount of time, prior to $t$, that the process spends in $j$.
From the large deviation estimates for the Markov chains it easily follows that for $c' \in (0,1)$   there are $\delta(c') > 0$ and $\varepsilon_0 > 0$ such that
\begin{equation} \label{u2}
\mathrm{P}_i(t^\varepsilon_{[(1-c')n(\varepsilon)]} > t(\varepsilon)) < e^{-\delta(c') n(\varepsilon)},~~~\varepsilon \leq \varepsilon_0.
\end{equation}
\begin{equation} \label{u1}
\mathrm{P}_i(t^\varepsilon_{[(1+c')n(\varepsilon)]} < t(\varepsilon)) < e^{-\delta(c') n(\varepsilon)},~~~\varepsilon \leq \varepsilon_0,
\end{equation}
Moreover, for $c' \in (0, c \wedge 1)$, there are $\delta(c,c')> 0$ and $\varepsilon_0 > 0$ such that
\begin{equation} \label{u3}
\mathrm{P}_i(\beta(j, t^\varepsilon_{[(1-c')n(\varepsilon)]}) < (1-c) \mu(j, \varepsilon) t(\varepsilon)) < e^{-\delta(c, c') n(\varepsilon)},~~~\varepsilon \leq \varepsilon_0,
\end{equation}
while for $c'> 0$ satisfying $(1+c)(1-c')/(1+c') > 1$, there are $\delta(c,c')> 0$ and $\varepsilon_0 > 0$ such that
\begin{equation} \label{u4}
\mathrm{P}_i(\beta(j, t^\varepsilon_{[(1+c')n(\varepsilon)]}) > (1+c) \mu(j, \varepsilon) t(\varepsilon)) < e^{-\delta(c, c') n(\varepsilon)},~~~\varepsilon \leq \varepsilon_0,
\end{equation}
We obtain (\ref{w3}) by combining (\ref{u2}) and (\ref{u3}). We obtain (\ref{w3b}) by combining (\ref{u1}) and (\ref{u4}).

From (\ref{w3b}) and the strong Markov property of the process, it follows that
\[
{\mathrm{P}_i}(\alpha(j, t(\varepsilon)) > k (1+c)\mu(j, \varepsilon) ) \leq e^{-k \delta(c) t(\varepsilon)/\bar{T}(\varepsilon)},~~~\varepsilon \leq \varepsilon_0,
\]
for each $k \in \mathbb{N}$. Combined with (\ref{w3}), this immediately implies (\ref{w2}).
\qed
\\

Next, let us consider the behavior of an asymptotically regular chain that is stopped when it enters a non-empty set $E \subseteq S$. Let
$\sigma = \inf\{t: X^\varepsilon_t \in E \}$, $\tau = \min\{n: Z_n \in E \}$, $\tau' = \min\{n: Z^\varepsilon_n \in E \}$.
\begin{lemma}  \label{jjj} Let $X^\varepsilon_t$ be asymptotically regular and $N \geq 2$. Suppose that the skeleton chain has one ergodic class and no transient states. Let $E$ be a non-empty subset of $S$.
Then for each $i \in S$ and $j \in E$,
\begin{equation} \label{tauone}
\lim_{\varepsilon \downarrow 0}  \mathrm{P}_i  ({X}^\varepsilon_\sigma = j ) =   \mathrm{P}_i  ({Z}_\tau = j ).
\end{equation}
If $t(\varepsilon)$ is such that $t(\varepsilon) \ll T(i, \varepsilon)$ for $i \in E$, $T(i, \varepsilon) \ll t(\varepsilon)$ for $i \notin E$, then
\begin{equation} \label{tautwo}
\lim_{\varepsilon \downarrow 0}  \mathrm{P}_i  ({X}^\varepsilon_{t(\varepsilon)} = j ) =   \mathrm{P}_i  ({Z}_\tau = j ).
\end{equation}
\end{lemma}
\proof Recall that $Z^\varepsilon_n$ is the discrete-time Markov chain obtained from $X^\varepsilon_t$ by discretizing time. Then
\[
\mathrm{P}_i  ({X}^\varepsilon_\sigma = j) = \mathrm{P}_i  ({Z}^\varepsilon_{\tau'} = j ) \rightarrow \mathrm{P}_i  ({Z}_\tau = j)~~{\rm as}~\varepsilon \downarrow 0,
\]
where the convergence follows from (\ref{prox}). Now let us prove (\ref{tautwo}). Given $\delta > 0$, find $k$ such that $\mathrm{P}_i(\tau > k) < \delta$. From the convergence of $Z^\varepsilon_n$ to $Z_n$, it then follows that
\[
\mathrm{P}_i(\tau' > k) < \delta
\]
for all sufficiently small $\varepsilon$. From the condition $T(i, \varepsilon) \ll t(\varepsilon)$ for $i \notin E$, it follows that
\[
\mathrm{P}_i(\sigma > t(\varepsilon)) < 2 \delta
\]
for all sufficiently small $\varepsilon$. Notice also that for $i \in E$,
\[
\mathrm{P}_i(t_1^\varepsilon \leq t(\varepsilon)) < \delta
\]
for all sufficiently small $\varepsilon$, as follows from the condition that $t(\varepsilon) \ll T(i, \varepsilon)$ for $i \in E$.
Since $\delta$ was arbitrary, using the last two inequalities, (\ref{tauone}), and the strong Markov property of the process, we obtain (\ref{tautwo}).
\qed
\\
\\
{\bf Remark.} The quantity $\mathrm{P}_i  ({X}^\varepsilon_\sigma = j )$ can be represented in terms of $i$-graphs (see Chapter 6 of \cite{FW}). Such a representation could be used as an alternative way to prove Lemma~\ref{jjj}.

\section{Metastable distributions for completely asymptotically regular families.} \label{se3}

\subsection{Formulation of the result.}

Suppose that $X^\varepsilon_t$ is completely asymptotically regular.
In this section, we'll study the distribution of  $X^\varepsilon_{t(\varepsilon)}$ at time scales $t(\varepsilon)$ that vary with $\varepsilon$ as $\varepsilon \downarrow 0$. The initial state $X^\varepsilon_0$ is assumed to be fixed.

To give a clearer picture, we first formulate the result in a particular case, with the general case to follow.
\begin{theorem} \label{maintheA}
Let $X^\varepsilon_t$ be asymptotically regular and $N \geq 2$. Suppose that the skeleton chain has one ergodic class and no transient states.
Suppose that for each  $1 \leq i \leq N$ either $t(\varepsilon) \ll T(i, \varepsilon)$ or $t(\varepsilon) \gg T(i, \varepsilon)$. Then there is a family of probability measures $\nu(i, \cdot)$, $i \in S$, on $S$ such that
\[
\lim_{\varepsilon \downarrow 0} {\mathrm{P}_i}(X^\varepsilon_{t(\varepsilon)} = j ) = \nu(i, j).
\]
\end{theorem}
The measure $\nu(i, \cdot)$ will be referred to as the metastable distribution for the initial state $i$ at the time scale $t(\varepsilon)$.
\proof Let $E = \{ i\in S: t(\varepsilon) \ll T(i, \varepsilon) \}$. If $E = \o$, then the result follows from (\ref{w1}), i.e., the limiting measure is the invariant measure, for each initial state. If $E \neq \o$, then the result follows from (\ref{tautwo}), i.e., the limiting measure is concentrated on $E$ and may depend on the initial point $i$.
\qed
\\


Now let formulate the result in the general case and describe how to identify the metastable distributions.
For $0 \leq r \leq \rho$,
the {\it inverse transition rate} of $S^{r}_i$ is defined as
\begin{equation} \label{itr}
T^{r}(i, \varepsilon) = ({\sum_{ j \neq i} Q^{r}_{ij}(\varepsilon)})^{-1},
\end{equation}
where the sum is over $1 \leq j \leq n_r$, $j \neq i$ (if $S^r_i \in S^{r+1}_k$ and $S^{r+1}_k \neq \{S^r_i\}$, this is asymptotically equivalent to taking only such $j$ that $S^{r}_j  \in S^{r+1}_k$). Note that $T^{r}(i, \varepsilon) \equiv +\infty$ if and only if $r = \rho$.  Our main result is the following.
\begin{theorem} \label{mainthe}
Suppose that for each $0 \leq r \leq \rho-1$ and each $1 \leq i \leq n_{r}$ either $t(\varepsilon) \ll T^{r}(i, \varepsilon)$ or $t(\varepsilon) \gg T^{r}(i, \varepsilon)$. Then there is a family of probability measures $\nu(i, \cdot)$, $i \in S$, on $S$ such that
\[
\lim_{\varepsilon \downarrow 0} {\mathrm{P}_i}(X^\varepsilon_{t(\varepsilon)} = j ) = \nu(i, j).
\]
\end{theorem}

The proof of the theorem will be based on two lemmas that we formulate next.
Let
\begin{equation} \label{invmy}
\mu^{r,k}(i) = \lim_{\varepsilon \downarrow 0} \mu^{r,k}(i, \varepsilon)~{\rm if}~S^r_i \in S^{r+1}_k;~~~\mu^{r,k}(i) = 0~{\rm if}~S^r_i \notin S^{r+1}_k,
\end{equation}
where the existence of the limit is guaranteed by Lemma~\ref{hjjh}.
 For $-1 \leq r < \rho$ and $j \prec S^{r+1}_k$ such that $j = S^0_{j_0} \in S^1_{j_1} ... \in S^{r+1}_{j_{r+1}} = S^{r+1}_k$, define
\begin{equation} \label{lmea}
\nu^{r} (j) = \mu^{0,j_1}(j_0) \mu^{1,j_2}(j_1)  ... \mu^{r,j_{r+1}}(j_{r}),
\end{equation}
where the right hand side is defined to be one if $r = -1$. For $i \in S$, let $r(i)$ be the minimal value of $r$ such that there is $k$ with $i \prec S^{r+1}_k$ and $t(\varepsilon) \ll T^{r+1}(k, \varepsilon)$.
Let $L(i) = \{l: S^{r(i)}_l \in S^{r(i)+1}_k~{\rm and}~t(\varepsilon) \ll T^{r(i)}(l, \varepsilon)\}$. Obviously, $L(i) = \o$ if
$r(i) = -1$ since then there are no $l$ for which $S^{r(i)}_l$ is defined.
The following lemma (proved in the next subsection) provides the description of the metastable distribution in some cases.
\begin{lemma} \label{le1a} Let the assumption made in Theorem~\ref{mainthe} hold.
Let $k$ be such that $i \prec S^{r(i)+1}_k$. Suppose that $L(i) = \o$. If $j$ satisfies (\ref{inclu}), then
\begin{equation} \label{uju}
\lim_{\varepsilon \downarrow 0} {\mathrm{P}_i}(X^\varepsilon_{t(\varepsilon)} = j ) = \nu^{r(i)} (j).
\end{equation}
If $j \prec S^{r(i)+1}_k$ does not hold, then $\lim_{\varepsilon \downarrow 0} {\mathrm{P}_i}(X^\varepsilon_{t(\varepsilon)} = j ) = 0$.
\end{lemma}

Now consider the case when $L = L(i) \neq \o$. Let $l(i)$ be such that $i \prec S^{r}_{l(i)}$, where $r = r(i)$. Observe that ${l(i)} \notin L(i)$.  Recall that the state space of $Y^{r, k, \varepsilon}_t$ is $S^{r+1}_k$. We define a new state space $\widetilde{S}^{r+1}_k$ by removing all the states
$S^{r}_l$, $l \in L$, from $S^{r+1}_k$ and adjoining the set $E = \{j: j \prec S^{r}_l~{\rm for}~{\rm some}~l \in L \}$. On this new state space, we
define the Markov chain $\widetilde{Y}^{\varepsilon}_t$. Its transition rates are defined as follows. If  $l, m \notin L$, $l \neq m$, are such that $S^{r}_l, S^{r}_m \in S^{r+1}_k$, then
\[
\widetilde{Q}_{S^{r}_l S^{r}_m}(\varepsilon) = Q^{r}_{lm}(\varepsilon).
\]
The transition rates between $S^{r}_l \subset S^{r+1}_k$, $l \notin L$, and $j \in E$ are defined as
\[
\widetilde{Q}_{S^{r}_l j}(\varepsilon) = \widetilde{Q}^{r}_{l j}(\varepsilon),
\]
where the quantity in the right hand side has been defined after the construction of the hierarchy.
The transition rate from $j \in E$ to any other state is zero, i.e., $E$ is the terminal set.
It is not difficult to show that $\widetilde{Y}^{\varepsilon}_t$ satisfies the assumptions of Lemma~\ref{jjj}, i.e., it coincides with an asymptotically regular Markov chain stopped upon entering $E$. (The proof of this statement is the same as the proof of Lemma~\ref{lereg0}.) Therefore, by (\ref{tauone}), there is a probability measure $\eta(i, \cdot)$ on $E$ such that
\[
\lim_{\varepsilon \downarrow 0} \lim_{t \rightarrow \infty} \mathrm{P}  (\widetilde{Y}^{\varepsilon}_t = j | \widetilde{Y}^{\varepsilon}_0 = S^{r}_{l(i)}) = \eta(i, j), ~~~j \in E.
\]
The proof of the following lemma is similar to that of Lemma~\ref{le1a}, and therefore not presented here.
\begin{lemma} \label{le1b} Let the assumption made in Theorem~\ref{mainthe} hold.
Suppose that  $j \prec S^{r(i)}_m$, where  $m \in L(i)$. Then
\[
\lim_{\varepsilon \downarrow 0} {\mathrm{P}_i}(X^\varepsilon_{t(\varepsilon)} = j ) = \sum_{i' \prec S^{r(i)}_m} \eta(i, i') \lim_{\varepsilon \downarrow 0} {\mathrm{P}_{i'}}(X^\varepsilon_{t(\varepsilon)} = j ),
\]
provided that the limits in the right hand side exist.
If $j$ is such that $j \prec S^{r(i)}_m$ does not hold for any $m \in L$, then
\[
\lim_{\varepsilon \downarrow 0} {\mathrm{P}_i}(X^\varepsilon_{t(\varepsilon)} = j ) = 0.
\]
\end{lemma}

\subsection{Proof of the main result.}

It is clear that Lemmas~\ref{le1a} and \ref{le1b} imply Theorem~\ref{mainthe}. Indeed, if $L(i) = \o$, then
the metastable distribution is given by Lemma~\ref{le1a}.
If $L(i) \neq \o$, then, by Lemma~\ref{le1b}, $\nu(i, j)$ is either equal to zero or is equal to a linear combination of the quantities
$ \nu(i', j)$ with $r(i') < r(i)$ (provided that all $\nu(i', j)$ are defined). The values of $\nu(i',j)$ can be found from Lemma~\ref{le1a} (when $L(i') = \o$), or again expressed in terms of metastable distributions with different initial points, using Lemma~\ref{le1b}. This recursive procedure can be continued until all the resulting initial points $j$ satisfy $L(j) = \o$, in which case Lemma~\ref{le1a} can be applied (which will happen in no more than $r(i)+1$ steps since $L(j) = \o$ whenever $r(j) = -1$).

It remains to prove Lemmas~\ref{le1a} and \ref{le1b}.
We give the proof of Lemmas~\ref{le1a} and omit the proof of Lemma~\ref{le1b}, since it is quite similar.
\\
\\
{\it Proof of Lemma~\ref{le1a}.} Let us start by briefly explaining the main idea of the proof. First, let us ``reduce'' the state space of $X^\varepsilon_t$ by clumping all the states $i$ with the property that $i \prec S^{r(i)}_m$ into a single state (recall the definition of '$\prec$' from (\ref{inclu})). The resulting process is well-approximated by the Markov chain $Y^{r(\varepsilon),k, \varepsilon}_t$, where $k$ is such that $X^\varepsilon_0 \prec S^{r(i)+1}_k$. Let us observe the process on a time scale $s(\varepsilon) \sim t(\varepsilon)$, such that $s(\varepsilon)$ is slightly smaller than $t(\varepsilon)$. Take $m$ such that $j \prec S^{r(i)}_m$. From the properties of $Y^{r(\varepsilon),k, \varepsilon}_t$ it then follows that $\lim_{\varepsilon \downarrow 0} {\mathrm{P}_i}(X^\varepsilon_{s(\varepsilon)} \prec S^{r(i)}_m ) = \mu^{r(i),k} (m)$, yielding the last factor in the expression (\ref{lmea}). Next, we can consider the process $X^\varepsilon_t$ on the time scale $t(\varepsilon) - s(\varepsilon)$ starting at a point $i' = X^\varepsilon_{s(\varepsilon)} \prec S^{r(i)}_m$. Provided that $s(\varepsilon)$ is chosen appropriately, the problem of identifying the limiting distribution at this time scale is similar to the original one, but
with $r(i') = r(i) -1$. Iterating the argument $r(i)$ times, we'll get the desired distribution. Let us now make the above arguments formal.

For $0 \leq r < \rho$, define the process $\bar{Y}^{r,\varepsilon}_t$ via
\[
\bar{Y}^{r,\varepsilon}_t = S^r_m~~~{\rm if}~~{X}^{\varepsilon}_t  \prec S^{r}_m.
\]
This is not necessarily a Markov process. However, below we'll see that on certain time scales it is close, in a certain sense, to the Markov process $Y^{r, k, \varepsilon}_t$, where $k$ is such that ${X}^{\varepsilon}_0 = i   \prec S^{r+1}_k$. Suppose that $s(\varepsilon)$ is such that
$s(\varepsilon) \ll T^{r+1}(k, \varepsilon)$, yet $T^r(l, \varepsilon) \ll s(\varepsilon)$ whenever $S^r_l \in S^{r+1}_k$. We claim (and will prove below) that
\begin{equation} \label{stee}
\lim_{\varepsilon \downarrow 0} \mathrm{P}(\bar{Y}^{r,\varepsilon}_{s(\varepsilon)} = S^r_l | X^\varepsilon_0 = i) = \mu^{r,k}(l)
\end{equation}
for each $1 \leq l \leq n_r$, where $\mu^{r,k}$, defined in (\ref{invmy}), is the limit, as $\varepsilon \downarrow 0$, for the invariant measures of the processes $Y^{r, k, \varepsilon}_t$.

Recall that $t(\varepsilon) \ll T^{r(i)+1}(k, \varepsilon)$, while $T^{r(i)}(l, \varepsilon) \ll t(\varepsilon)$ whenever $S^{r(i)}_l \in S^{r(i)+1}_k$. Choose $s(\varepsilon) < t(\varepsilon)$ with the following properties:

(a) $s(\varepsilon) \ll T^{r(i)+1}(k, \varepsilon)$, while $T^{r(i)}(l, \varepsilon) \ll s(\varepsilon)$ whenever $S^{r(i)}_l \in S^{r(i)+1}_k$.

(b) For each $l$ such that $\mu^{r(i), k}(l) > 0$, we have $t(\varepsilon) - s(\varepsilon) \ll T^{r(i)}(l, \varepsilon)$, while
$T^{r'}(l', \varepsilon) \ll t(\varepsilon) - s(\varepsilon)$ whenever $S^{r'}_{l'}$ is such that $S^{r'}_{l'} \prec S^{r(i)}_l$ and there is $j \prec S^{r(i)}_l$ for which $j \prec S^{r'}_{l'}$ does not hold. (Here $S^{r'}_{l'} \prec S^{r(i)}_l$ means that $S^{r'}_{l'} \in S^{r'+1}_{l_1} \in ... \in S^{r'+m}_{l_m} \in S^{r(i)}_l$ for some $m$, $l_1,...,l_m$.)

The existence of such $s(\varepsilon)$ follows from Lemma~\ref{hjjh}, which guarantees that the quantities $T^{r(i)}(l, \varepsilon)$ (with $l$ such that $\mu^{r(i), k}(l) > 0$) are all asymptotically equivalent, up to multiplicative constants.

By property (a), from (\ref{stee}) with $r = r(i)$, it follows that  ${X}^{\varepsilon}_{s(\varepsilon)}  \prec S^{r(i)}_l$ with probability close to $\mu^{r(i),k}(l)$ for $l$ such that $\mu^{r(i), k}(l) > 0$. Then, by the Markov property, ${X}^{\varepsilon}_{s(\varepsilon)}$ can be taken
as a new starting point for the process that will be studied on the time interval $t(\varepsilon) - s(\varepsilon)$. For this new time scale and the new initial point, all the assumptions of Lemma~\ref{le1a} are satisfied, with the exception that the value of $r(i)$ is reduced by at least one. Thus, conditioning on the event that $\bar{Y}^{r(i),\varepsilon}_{s(\varepsilon)} = S^r_l$, with $l$ such that $j \prec S^{r(i)}_l$, gives
\[
\lim_{\varepsilon \downarrow 0} {\mathrm{P}_i}(X^\varepsilon_{t(\varepsilon)} = j ) = \mu^{r(i),k}(l) \lim_{\varepsilon \downarrow 0} {\mathrm{P}_{i'}}(X^\varepsilon_{t(\varepsilon) - s(\varepsilon)} = j ),
\]
provided that the limit in the right hand side is defined and does not depend on $i' \prec S^{r(i)}_l$. Iterating this argument $r(i)+1$ times gives (\ref{uju}). It remains to prove (\ref{stee}).

For $i \in S$ and $-1 \leq r < \rho$, let $k$ be such that $i \prec S^{r+1}_k$. Let $A^{r} = \{j: j \prec S^{r+1}_k \}$. For $0 \leq r < \rho$, let $ \tilde{X}^{r ,\varepsilon}_t$ be the Markov chain on the state space $A^{r}$, whose
transition rates agree with those of the chain $X^\varepsilon_t$ (i.e., $ \tilde{X}^{r ,\varepsilon}_t$ is obtained from $X^\varepsilon_t$ by disallowing transitions outside the specified state space).

 For $j \prec S^{r+1}_k$ such that $j = S^0_{j_0} \in S^1_{j_1} ... \in S^{r+1}_{j_{r+1}}$ with $j = j_0$, $k = j_{r+1}$, let
\[
M^r(j, \varepsilon) = \mu^{0,j_1}(j_0, \varepsilon) \mu^{1,j_2}(j_1,\varepsilon)  ... \mu^{r,j_{r+1}}(j_{r},\varepsilon).
\]
For $0 \leq r < \rho$, let $\alpha^r(j, s(\varepsilon))$ be the proportion of time, prior to $s(\varepsilon)$, that the process $\tilde{X}^{r,\varepsilon}_{t}$ spends in $j$.
Let
\[
\bar{T}^r(\varepsilon) = \sum_{l: S^r_l \in S^{r+1}_k} \mu^{r,k}(l) T^r(l,\varepsilon).
\]
Define the process $ \tilde{Y}^{r,\varepsilon}_t$ via
\[
{\tilde{Y}}^{r,\varepsilon}_t = S^r_m~~~{\rm if}~~{\tilde{X}}^{r,\varepsilon}_t  \prec S^{r}_m.
\]
Recall that $s(\varepsilon)$ satisfies
$s(\varepsilon) \ll T^{r+1}(k, \varepsilon)$ and $T^r(l, \varepsilon) \ll s(\varepsilon)$ whenever $S^r_l \in S^{r+1}_k$.
 We claim that for each $0 \leq r < \rho$ and $s(\varepsilon)$ satisfying these conditions, the following asymptotic relations hold. For each $i, j \in A^r$ and $l$ such that $S^r_l \in S^{r+1}_k$,
 \begin{equation} \label{steexx}
\lim_{\varepsilon \downarrow 0} \mathrm{P}({\tilde{Y}}^{r,\varepsilon}_{s(\varepsilon)} = S^r_l | \tilde{X}^{r,\varepsilon}_0 = i) = \mu^{r,k}(l).
\end{equation}
For each $c > 0$, there are $\delta(c) > 0$ and a function $\varphi(\varepsilon) > 0$ such that $\lim_{\varepsilon \downarrow 0} \varphi(\varepsilon) = + \infty$ and
\begin{equation} \label{w3xx}
{\mathrm{P}_i}(\alpha^r(j, s(\varepsilon)) < (1-c)M^r(j, \varepsilon) ) \leq e^{-\delta(c) \varphi(\varepsilon)},
\end{equation}
\begin{equation} \label{w3bxx}
{\mathrm{P}_i}(\alpha^r(j, s(\varepsilon)) > (1+c)M^r(j, \varepsilon) ) \leq e^{-\delta(c) \varphi(\varepsilon)}.
\end{equation}
When $r =0$, (\ref{steexx})-(\ref{w3bxx}) follow from Lemma~\ref{interm}. Next, let us sketch the inductive step, i.e., the proof that relations
(\ref{steexx})-(\ref{w3bxx}) hold for a fixed ${r} > 0$, assuming that they hold for all the smaller values of $r$.

The process ${\tilde{Y}}^{r,\varepsilon}_t$ is close to the Markov process $Y^{r, k, \varepsilon}_t$ in the following sense.

(a) If $l \neq m$ and $S^{r}_l, S^{r}_m \in S^{r+1}_k$, $i \prec S^{r}_l$, then
\[
\mathrm{P}({\rm the}~{\rm first}~{\rm jump}~{\rm of}~{\tilde{Y}}^{r,\varepsilon}_t~{\rm is}~{\rm to}~S^r_m|\tilde{X}^{r ,\varepsilon}_0 = i) \sim Q^{r}_{lm}(\varepsilon)~~{\rm as}~\varepsilon \downarrow 0.
\]

(b) If $\beta^\varepsilon$ is the random time till the first transition of ${\tilde{Y}}^{r,\varepsilon}_t$, then
\[
{\mathrm{E}} (\beta^\varepsilon| \tilde{X}^{r ,\varepsilon}_0 = i) \sim  T^{r}(l, \varepsilon),~~~{\rm as}~\varepsilon \downarrow 0,
\]
and there are $\delta > 0$ and $\varepsilon_0 > 0$ such that
\[
{\mathrm{P}}(\beta^\varepsilon \geq \lambda  T^{r}(l, \varepsilon)|\tilde{X}^{r ,\varepsilon}_0 = i) \leq  e^{-\delta  \lambda},~~~~\lambda \geq 2, ~~\varepsilon \in (0, \varepsilon_0].
\]
The validity of (a) and (b) easily follows by examining the Markov chain ${\tilde{X}}^{r-1,\varepsilon}_t$ and utilizing the fact that (\ref{w3xx})-(\ref{w3bxx}) hold with $r$ replaced by $r-1$.

Lemma~\ref{interm} (with $s(\varepsilon)$ instead of $t(\varepsilon)$) can be applied to the process $Y^{r, k, \varepsilon}_t$. However, we are interested a similar result for the process ${\tilde{Y}}^{r,\varepsilon}_t$. It is not difficult to
see that conditions (a) and (b) on the transition probabilities and transition times are sufficient for the proof of Lemma~\ref{interm} to go through and for the result to be valid for the process ${\tilde{Y}}^{r,\varepsilon}_t$. Thus we have (\ref{steexx}).
Moreover, for $i \prec S^r_l$ and each $c > 0$, there are $\delta(c) > 0$ and a function $\varphi(\varepsilon) > 0$ such that $\lim_{\varepsilon \downarrow 0} \varphi(\varepsilon) = + \infty$ and
\[
{\mathrm{P}}({\tilde{\alpha}}^r(l, s(\varepsilon)) < (1-c)\mu^{r,k}(l, \varepsilon) | \tilde{X}^{r ,\varepsilon}_0 = i) \leq e^{-\delta(c) \varphi(\varepsilon)},
\]
\[
{\mathrm{P}}({\tilde{\alpha}}^r(l, s(\varepsilon)) > (1+c)\mu^{r,k}(l, \varepsilon) | \tilde{X}^{r ,\varepsilon}_0 = i) \leq e^{-\delta(c) \varphi(\varepsilon)},
\]
where ${\tilde{\alpha}}^r(l, s(\varepsilon))$ is the proportion of time, prior to $s(\varepsilon)$, that the process ${\tilde{Y}}^{r,\varepsilon}_t$
spends in the state $S^r_l$. Together with
(\ref{w3xx})-(\ref{w3bxx}) for $r-1$ instead of $r$, these are easily seen to control the proportion of time, prior to $s(\varepsilon)$, that $\tilde{X}^{r,\varepsilon}_{t}$ spends in $j$, thus yielding (\ref{w3xx})-(\ref{w3bxx}) for $r$.

From (\ref{w3xx})-(\ref{w3bxx}) it easily follows that
\[
\mathrm{P}_i(X^\varepsilon_t \prec S^{r+1}_k~~{\rm for}~{\rm all}~t \leq s(\varepsilon)) \rightarrow 1~~~{\rm as}~\varepsilon \downarrow 0.
\]
Thus (\ref{stee}) follows from (\ref{steexx}). \qed

\section{Complete asymptotic regularity} \label{car}
In this section, we prove Lemma~\ref{lereg0} and briefly discuss a couple of generalizations. We start with the following simple lemma.
\begin{lemma} \label{slem}
 Suppose that $a_1(\varepsilon),...,a_n(\varepsilon)$ and $b_1(\varepsilon),...,b_{n'}(\varepsilon)$ are positive functions such that there are limits
\[
\lim_{\varepsilon \downarrow 0} \frac{a_i(\varepsilon)}{b_{i'}(\varepsilon)} \in [0,\infty],~~1 \leq i \leq n,~1 \leq i' \leq n'.
\]
Let $A(\varepsilon) = a_1(\varepsilon) + ... + a_n(\varepsilon)$, $B(\varepsilon) = b_1(\varepsilon) + ... + b_{n'}(\varepsilon)$.
Then there is the limit
\[
\lim_{\varepsilon \downarrow 0} \frac{A(\varepsilon)}{B(\varepsilon)}\in [0, \infty].
\]
\end{lemma}
\proof
For each $i'$, there is the limit
\[
\lim_{\varepsilon \downarrow 0} \frac{A(\varepsilon)}{b_{i'}(\varepsilon)}  = \lim_{\varepsilon \downarrow 0} \frac{a_1(\varepsilon)}{b_{i'}(\varepsilon)} + ... +
\lim_{\varepsilon \downarrow 0}  \frac{a_n(\varepsilon)}{b_{i'}(\varepsilon)} \in [0, \infty].
\]
Therefore, there is the limit
\[
\lim_{\varepsilon \downarrow 0} \frac{A(\varepsilon)}{B(\varepsilon)} = (\lim_{\varepsilon \downarrow 0} \frac{b_1(\varepsilon)}{A(\varepsilon)} + ... +\lim_{\varepsilon \downarrow 0}  \frac{ b_{n'} (\varepsilon)}{A(\varepsilon)} )^{-1} \in [0, \infty].
\]
\qed
\\
\\
\\
{\it Proof of Lemma~\ref{lereg0}.}
The case $N = 1$ is trivial. Let us assume that $N \geq 2$. Recall the decomposition (\ref{decomp}) of $S$ into ergodic classes and transient states.
To prove the lemma, we need to show that for $a \geq 1$ there is the limit
\begin{equation} \label{lpro}
\lim_{\varepsilon \downarrow 0} \left( \frac{Q_{k_1 l_1}(\varepsilon)}{Q_{m_1 n_1}(\varepsilon)} \times ... \times \frac{Q_{k_a l_a}(\varepsilon)}{Q_{m_a n_a}(\varepsilon)}  \right) \in [0, \infty],
\end{equation}
provided that $k_1 \neq l_1,...,k_a \neq l_a, m_1 \neq n_1,...,m_a \neq n_a$.. We will repeatedly use Lemma~\ref{slem}, which will allow us to replace each of the factors above by simpler expressions. First consider a factor of the form $Q_{kl}/Q_{mn}$ under the assumption that $S_k = \{i\}$ and $S_m = \{i'\}$, i.e., $S_k$ and $S_m$ have only one element each. In this case
\[
\frac{Q_{kl}(\varepsilon)}{Q_{mn}(\varepsilon)} = \frac{\sum_{j \in S_l} q_{ij}(\varepsilon)  }{\sum_{j' \in S_n} q_{i'j'}(\varepsilon)}.
\]
Thus, by Lemma~\ref{slem}, it is sufficient to prove the existence of the limit in (\ref{lpro}) with all such factors replaced by those of the form $q_{ij}/q_{i'j'}$.

Next consider a factor of the form $Q_{kl}/Q_{mn}$ under the assumption that one of the sets $S_k$ and $S_m$ (say, $S_k$) has at least two elements, while the other one has one element, i.e., $S_m = \{i'\}$. Then the chain
$Y^{k,\varepsilon}_t$ satisfies the assumptions of Lemma~\ref{hjjh}, and therefore, by (\ref{defq}),
\[
\frac{Q_{kl}(\varepsilon)}{Q_{mn}(\varepsilon)} = \frac{\sum_{i \in S_k} \sum_{j \in S_l} \mu^{k}(i, \varepsilon) q_{ij}(\varepsilon)}{ \sum_{j' \in S_n} q_{i'j'}(\varepsilon)} \sim
\frac{\sum_{i \in S_k} \sum_{j \in S_l} \lambda(i) T(i, \varepsilon) q_{ij}(\varepsilon) }{ \sum_{j' \in S_n} q_{i'j'}(\varepsilon)  \sum_{i'' \in S_k}( \lambda(i'') T(i'',\varepsilon))},
\]
where $\lambda$ and $T$ are the invariant measure for the skeleton chain and the inverse transition rate, respectively, for the chain
$Y^{k,\varepsilon}_t$. Thus,  by Lemma~\ref{slem}, it is sufficient to prove the existence of the limit in (\ref{lpro}) with $Q_{kl}/Q_{mn}$ replaced by $T(i, \varepsilon) q_{ij}(\varepsilon)(T(i'', \varepsilon) q_{i'j'}(\varepsilon))^{-1}$, where $i, i'' \in S_k$. By the definition of $T$,
\[
\frac{T(i, \varepsilon) q_{ij}(\varepsilon)}{T(i'', \varepsilon) q_{i'j'}(\varepsilon)} =
\frac{{\sum_{b' \in S_k, b' \neq i''} q_{i''b'}(\varepsilon)} q_{ij}(\varepsilon)}{{\sum_{b \in S_k, b \neq i} q_{ib}(\varepsilon)} q_{i'j'}(\varepsilon)}.
\]
By Lemma~\ref{slem}, each such expression can be replaced by $( q_{i''b'}(\varepsilon) q_{ij}(\varepsilon))/ ( q_{ib}(\varepsilon) q_{i'j'}(\varepsilon))$.

The final case, when $S_k$ and $S_m$ have at least two elements each, is treated similarly, resulting in $Q_{kl}/Q_{mn}$ being replaced by a product of three factors of the form $q_{ij}/q_{i'j'}$.
Thus we see that each factor in (\ref{lpro}) can be replaced by either one, two, or three factors of the form $q_{ij}/q_{i'j'}$.  Therefore, the limit in (\ref{lpro}) exists since the original chain is completely asymptotically regular. \qed
\\


\section{Remarks and generalizations.} \label{gene}

{\bf (A)} One could replace the complete asymptotic regularity by a somewhat weaker assumption that also implies the asymptotic regularity of all the reduced chains appearing in the inductive construction of the hierarchy.
We will say that an asymptotically regular family of Markov chains $X^\varepsilon_t$ satisfies Condition (A) if
for each $0 \leq a \leq N-2$
the following finite or infinite limit exists
\[
\lim_{\varepsilon \downarrow 0} \left( \frac{q_{k_1 l_1}(\varepsilon)}{q_{k_1 n_1}(\varepsilon)} \times ... \times \frac{q_{k_a l_a}(\varepsilon)}{q_{k_a n_a}(\varepsilon)} \times \frac{q_{k' l'}(\varepsilon)}{q_{m' n'}(\varepsilon)} \right) \in [0, \infty],
\]
provided that $k_1,...,k_a$ are all distinct, $k', m' \notin \{k_1,...,k_a\}$, and $k_1 \neq l_1,..., k_a \neq l_a$, $k_1 \neq n_1,...,k_a \neq n_a$,
$k' \neq l', m' \neq n'$. The proof of the following lemma is similar to that of Lemma~\ref{lereg0}, and so we don't provide it here.
\begin{lemma} \label{lereg}
If $X^{\varepsilon}_t$  satisfies Condition (A), then the reduced chain also satisfies Condition~(A).
\end{lemma}
If $X^{\varepsilon}_t$ satisfies Condition~(A),
then, by Lemma~\ref{lereg}, so do the reduced Markov chains that appear at each step of the inductive construction of the hierarchy, which implies that all of them are asymptotically regular.
\\
\\
{\bf (B)} Next, let us mention that all the above analysis can be easily adapted to the case of discrete-time Markov chains. Instead of (\ref{contti}), we assume that $X^\varepsilon_n$, $n \geq 0$, is a Markov chain with transition probabilities $p_{ij}(\varepsilon)$, $i,j \in S$. We impose an additional assumption that $p_{ii}(\varepsilon) \geq c > 0$ for all $i \in S$,  $\varepsilon > 0$. This is needed in order for the discrete time analogue of Lemma~\ref{interm} to remain valid.

The definitions of asymptotic regularity and complete asymptotic regularity remain the same as in the continuous time case, with transition rates
$q_{ij}(\varepsilon)$ replaced by transition probabilities $p_{ij}(\varepsilon)$. The reduced Markov chains the chains
chains $Y^{r, k, \varepsilon}_t$ can be still defined in continuous time, simply replacing $q_{ij}(\varepsilon)$ by $p_{ij}(\varepsilon)$ in all
the definitions. The definition of the inverse transition rates (\ref{itr}) remains the same. The theorem on metastable distributions  now takes the following form.
\begin{theorem} \label{mainthe2}
Let $n: (0,\infty) \rightarrow \mathbb{N}$ be such that
for each $0 \leq r \leq \rho-1$ and each $1 \leq i \leq n_{r}$ either $n(\varepsilon) \ll T^{r}(i, \varepsilon)$ or $n(\varepsilon) \gg T^{r}(i, \varepsilon)$. Then there is a family of probability measures $\nu(i, \cdot)$, $i \in S$, on $S$ such that
\[
\lim_{\varepsilon \downarrow 0} {\mathrm{P}_i}(X^\varepsilon_{n(\varepsilon)} = j ) = \nu(i, j).
\]
\end{theorem}
\noindent
\\
The proof of this theorem is identical to that of Theorem~\ref{mainthe}.
\\
\\
{\bf (C)} Finally, consider an example of a completely asymptotically regular family of Markov chains. Given numbers $\alpha_{ij}$, $\beta_{ij}$, and $\gamma_{ij}$,
$i \neq j$, such that $\alpha_{ij} > 0$, assume that the transition rates satisfy
\begin{equation} \label{e1}
q_{ij}(\varepsilon) \sim \alpha_{ij} \varepsilon^{\beta_{ij}} \exp({-\gamma_{ij} \varepsilon^{-1}}),~~{\rm as}~\varepsilon \downarrow 0,~~~i \neq j.
\end{equation}
It is clear that these functions satisfy (\ref{care}), and therefore the corresponding family is completely asymptotically regular. Replacing $\varepsilon$ by $|\ln \tilde{\varepsilon}|^{-1}$, we obtain functions
\begin{equation} \label{e2}
\tilde{q}_{ij}(\tilde{\varepsilon}) \sim \alpha_{ij} |\ln \tilde{\varepsilon}|^{-\beta_{ij}}  (\tilde{\varepsilon})^{\gamma_{ij}} ,~~{\rm as}~\varepsilon \downarrow 0,~~~i \neq j,
\end{equation}
which also satisfy (\ref{care}).
The systems discussed in the Introduction lead to Markov chains with transition rates that satisfy either (\ref{e1}) or (\ref{e2}), with the exception that the condition $\alpha_{ij} > 0$ may be violated, i.e., some of the coefficients may be equal to zero. In fact, this positivity condition (or  condition (a) in the definition of asymptotic regularity) are not that crucial. If $i$ and $j$ are
such that $\alpha_{ij} = 0$, the transition rates can be re-defined  for those $(i,j)$ by taking $\alpha_{ij} = 1$ and $\gamma_{ij}$ sufficiently large, resulting in a completely asymptotically regular family with the same metastable behavior as the original one.
\\
\\

\noindent {\bf \large Acknowledgements}: While working on this
article, M. Freidlin was supported by NSF grant DMS-1411866
and L. Koralov was supported by NSF grant DMS-1309084.
\\
\\

\end{document}